% Fichier DehTeXSource.
%
% Ce fichier contient les macros du format DehTeX, ainsi qu'un commentaire
% qui peut etre edite seul et constituer une notice.
% Pour le compiler le fichier, enlever % de devant la ligne \compiler et placer
% % devant la ligne \editer; pour l'editer, faire l'inverse. Il faut compiler
% avant d'editer (car l'edition fait appel aux macros), et compiler avec
% la derniere version de DehTeX.
% Pour rajouter une macro \xx, taper sa definition a sa place 
% alphabetique, puis ajouter le commentaire par \comm(xx)blabla!!! (trois
% points d'exclamation).
% En entrant le commentaire, on tape \\ pour faire imprimer \, \# pour faire
% imprimer #, \{ et \} pour faire imprimer { et }.

\long\def\comm(#1)#2!!!{}

{\gdef\editer
                                {\tolerance 3000
                                \gdef\\{{\tt\char'134}}\gdef\#{{\tt\char'043}}          
        
                                \gdef\{{{\tt\char'173}}\gdef\}{{\tt\char'175}}
                                \long\gdef\comm(##1)##2!!!
                                                        
        {{\leftskip20mm\parindent=0mm\medskip\hskip-20mm
                                                                \setbox1=\hbox{\\\tt##1}
                                                                {\ifdim\wd1<17mm\hbox to 
19mm{\box1\hfill}%
                                                                \else\hbox to \wd1{\unhbox1} : 
\fi}##2.\par}}
                                \centerline{\bfXII Les commandes de DehTeX}\bigskip
                                \datef\bigskip\bigskip
                                Num\'eros des familles de fontes: 
        0=cmr$\fl$rm; 1=cmmi$\fl$mi; 2=cmsy$\fl$sy; 3=cmex;  4=cmti$\fl$it;
                                5=cmsl$\fl$sl; 6=cmbx$\fl$bf; 7=cmtt$\fl$tt; 8=cmb$\fl$be; 
                                9=cmmib$\fl$bm; A=msam$\fl$sa; B=msbm$\fl$sb; 
C=eufm$\fl$go;
                                D=eurm$\fl$eu; E non utilis\'e; F non utilis\'e;
                                \def\dump{}}}

%\compiler
%\editer

\def\<<{{\cyr\char'074~}}
                                \comm(<<)Ecrit \<<(avec blanc inseccable apr\`es)!!!

\def\>>{~{\cyr\char'076}}
                                \comm(>>)Ecrit \>> (avec blanc inseccable avant)!!!

\def\\{\hbox{$\backslash$}}
                                \comm($\\$)Ecrit $\\$; ne n\'ecessite pas de blanc apr\`es!!! 

\let\(=\langle
                                \comm(()Ecrit $\($; mode math\'ematique!!!

\let\)=\rangle
                                \comm({)})Ecrit $\)$; mode math\'ematique!!!

\def\0{{\mi 0}}\def\1{{\mi 1}}\def\2{{\mi 2}}\def\3{{\mi 3}}
\def\4{{\mi 4}}\def\5{{\mi 5}}\def\6{{\mi 6}}\def\7{{\mi 7}}
\def\8{{\mi 8}}\def\9{{\mi 9}}
                                \comm({\it chiffre})Ecrit le chiffre correspondant en 
\<<oldstyle\>>: \0,
                                                                \1, \2, \3, \4, \5, \6, \7, \8, \9!!!

\font\V=cmr5 \font\VI=cmr6 \font\VII=cmr7 \font\VIII=cmr8
\font\IX=cmr9 \font\XII=cmr10 scaled\magstep1
\font\XIV=cmr10 scaled\magstep2
\font\XVIII=cmr10 scaled\magstep3
\font\XXIV=cmr10 scaled\magstep4
                                \comm({{\it chiffre romain}})Appelle la taille de caract\`ere
                                correspondante dans la fonte cmr; disponibles: {\tt\\V, \\VI, \\VII,
                                \\VIII, \\IX, \\XII, \\XIV, \\XVIII, \\XXIV}!!!

\mathchardef\a="710B
                                \comm(a)Ecrit $\a$; mode math\'ematique!!!

\def\adots{\mathinner{\mkern2mu\raise1pt\hbox{.}
                                                                \mkern3mu\raise4pt\hbox{.}
                                                                \mkern1mu\raise7pt\hbox{.}}}
                                \comm(adots)Ecrit $\adots$; mode math\'ematique!!!

\let\al=\aleph
                                \comm(al)Ecrit $\al$; mode math\'ematique!!!

\def\and{\hbox{ and }}
                                \comm(and)Ecrit $\and$ en romain avec blanc avant et apr\`es, 
quel
                                                                que soit le mode (le m\^eme en 
pench\'e s'appelle {\tt\\sland})!!!

\def\ape{\mathchar"3A44}
                                \comm(ape)Ecrit $\ape$; mode math\'ematique!!!
                 
\def\app{\mathord{\in}}
                                \comm(app)Ecrit $\app$ sans blanc avant ni apr\`es; mode 
math\'ematique!!!

\let\aps=\triangleright
                                \comm(aps)Ecrit $\aps$; mode math\'ematique!!!
                
\let\av=\triangleleft
                                \comm(av)Ecrit $\av$; mode math\'ematique; autre nom: {\tt 
avs}!!!
                
\def\ave{\mathchar"3A45}
                                \comm(ave)Ecrit $\ave$; mode math\'ematique!!!

                                \comm(avs)Voir {\tt ave}!!!

\mathchardef\b="710C       
                                \comm(b)Ecrit $\b$; mode math\'ematique!!!

\def\bbar#1{\overline{\mathstrut#1}}
                                \comm(bbar)La commande {\tt\\bbar\#} place une barre sur \#, \`a
                                                                une hauteur uniforme et assez 
grande; mode math\'ematique!!!

%famille8: be;
                                                                \font\cmbX=cmb10 at 10pt 
\font\cmbVII=cmb10 at 7pt 
                                                                \font\cmbV=cmb10 at 5pt 
\textfont8=\cmbX
                                                        
        \scriptfont8=\cmbVII\scriptscriptfont8=\cmbV
                                                                \def\be{\fam8\cmbX}
                                                                \let\bfe=\be
                                \comm(be)Passe \`a la famille \<<gras \'etroit\>> (fontes cmb); 
{\be
                                                                Ceci est du be}; autre nom: 
{\tt\\bfe}!!!

%famille6: bf;
                                                        
        \font\cmbxX=cmbx10\font\cmbxVII=cmbx7 
                                                                \font\cmbxV=cmbx5 
\textfont6=\cmbxX
                                                        
        \scriptfont6=\cmbxVII\scriptscriptfont6=\cmbxV
                                                                \def\bf{\fam6\cmbxX}

\font\bfXII=cmbx10 scaled\magstep1  
              \font\bfXIV=cmbx10 scaled\magstep2  
              \font\bfXVIII=cmbx10 scaled\magstep3  
              \font\bfXXIV=cmbx10 scaled\magstep4
                                \comm(bf{{\it chiffre romain}})Passe \`a la fonte \<<gras\>>\ 
(cmbx)
                                                                de la taille correspondante;
                                                                disponibles: {\tt\\bfXII, \\bfXIV, 
\\bfXVIII, \\bfXXIV}!!!

\def\bg{\medbreak\medskip}
        \comm(bg)Provoque un grand saut et favorise la coupure; la moiti\'e du
                                saut est supprim\'ee apr\`es un display!!!

\def\bgni{\medbreak\medskip\noindent}
                                \comm(bgni)Provoque un grand saut, favorise la coupure et
                                supprime l'indentation apr\`es!!!

\let\bi=\cmbxtiX
                                \comm(bi)Passe \`a la fonte \<<gras italique\>> (fontes cmbxti); 
{\bi
                                Ceci est du bi}!!!

                                \comm(bl)Ecrit $\bullet$!!!

%famille9: bm;
                                                        
        \font\cmmibX=cmmib10\font\cmmibVII=cmmib7
                \font\cmmibV=cmmib5 \textfont9=\cmmibX
                \scriptfont9=\cmmibVII\scriptscriptfont9=\cmmibV
                \def\bm{\fam9\cmmibX}
                                \comm(bm)Appelle la famille \<<gras math\'ematique\>> (fontes
                                                                cmmib); {\bm Ceci est du bm}!!!

\def\bn{\bigskip\nobreak}
                                \comm(bn)Provoque un saut et d\'efavorise la coupure!!!

                                \comm(bnni)Provoque un saut, d\'efavorise la coupure et 
supprime
                                                                l'indentation!!!

\long\def\boxit#1#2{\setbox1=\hbox{\kern#1{#2}\kern#1}%
                                                                \dimen1=\ht1 \advance\dimen1 by 
#1
                                                                \dimen2=\dp1 \advance\dimen2 by 
#1
                                                                \setbox1=\hbox{\vrule 
height\dimen1 depth\dimen2\box1\vrule}%
                                                        
        \setbox1=\vbox{\hrule\box1\hrule}%
                                                                \advance\dimen1 by .4truept 
\ht1=\dimen1
                                                                \advance\dimen2 by .4truept 
\dp1=\dimen2 \box1\relax}
                                \comm(boxit)Encadre son second argument; syntaxe
                                                                {\tt\\boxit$\{$2pt$\}\{$xx$\}$} 
donne \boxit{2pt}{{\tt xx}}!!!

\let\bs=\cmbxslX
                                \comm(bs)Passe \`a la fonte \<<gras pench\'e\>> (fontes cmbxsl); 
{\bs
                                ceci est du bs}!!!

\def\card{{\rm card}}
                                \comm(card)Ecrit \card!!!               

\def\carre{\hbox{\sa\char'004}}
                                \comm(carre)Ecrit \carre!!!

\catcode`\@=11\def\cases#1{\left\{\,\vcenter{\normalbaselines\m@th
                                                        
        \ialign{$##\hfil$&\quad{\rm##}\hfil\crcr#1\crcr}}\right.}
                                                                \catcode`\@=12

\def\cf{\hbox{\it c.f. }}
                                \comm(cf)Ecrit \cf\ avec un espace apr\`es!!!

\font\cmdunhX=cmdunh10
                                \comm(cmdunhX)Passe \`a la fonte cmdunh10; {\cmdunhX Ceci 
est du
                                                                cmdunh10}!!!

\def\comp{{\scriptscriptstyle\circ}}
                                \comm(comp)Ecrit $\comp$; mode math\'ematique!!!

                                \comm(Cor)La commande {\tt Cor\#} provoque un saut, puis 
\'ecrit
                                                                sans indentation {\bf Corollary\#.-
}!!!

% \font\cyr=wncyr10
                                \comm(cyr)Fait passer dans la fonte \<<wncyr10\>>; {\cyr Ceci 
est du
                                                                cyr}!!!

\mathchardef\d="710E
                                \comm(d)Ecrit $\d$; mode math\'ematique!!!

\mathchardef\D="7101
                                \comm(D)Ecrit $\D$; mode math\'ematique!!!

\def\date{\par\line{\hfil\ifcase\month\or January\or
                                                                February\or March\or April\or 
May\or June\or July\or August\or
                                                                September\or October\or 
November\or December\fi\ 
                                                                {\oldstyle\the\day},\ 
{\oldstyle\the\year}}}
                                \comm(date)Met la date en anglais \`a droite de la ligne!!!

\def\datef{\par\line{\hfil le\ {\oldstyle\the\day}\
                 \ifcase\month\or
                 janvier\or f\'evrier\or mars\or avril\or mai\or juin\or
                 juillet\or ao\^ut\or septembre\or octobre\or novembre\or
                 d\'ecembre\fi\ {\oldstyle\the\year}}}
                                \comm(datef)Met la date en francais \`a droite de la ligne!!!

\def\Def{\bgni{\bf Definition.- }}
                                \comm(Def)Provoque un saut, puis \'ecrit {\bf Definition.-}!!!

\let\der=\partial
                                \comm(der)Ecrit $\der$; mode math\'ematique!!!

                                \comm(diagram)Pour coder un diagramme; m\^eme syntaxe que
                                                                {\tt\\matrix}, simplement les espaces 
sont red\'efinis!!!

\def\Dom{\mathord{\rm Dom}}
                                \comm(Dom)Ecrit $\Dom$; mode math\'ematique!!!

\mathchardef\e="7122
                                \comm(e)Ecrit $\e$; mode math\'ematique!!!

\def\eg{\hbox{\it e.g. }}
                                \comm(eg)Ecrit \eg avec un espace apr\`es!!!

                                \comm(entete)Place l'entete de Caen en haut de page; attention! il 
faut
                                                                que le sceau soit dans le fichier 
pictures!!!

\let\equ=\Longleftrightarrow
                                \comm(equ)Ecrit $\equ$; mode math\'ematique!!!

                                \comm(esp)Force un blanc, m\^eme dans les fontes o\`u l'espace 
est nul,
                                                                comme cmmi qui sert dans 
oldstyle!!!

\def\et{\hbox{ et }}
                                \comm(et)Ecrit $\et$ en tout mode (avec un espace avant et 
apr\`es)!!!

\def\etc{\hbox{\it etc\dots\ }}
                                \comm(etc)Ecrit \etc avec un espace apr\`es!!!

%famille 13: eu
                                                        
        \font\euX=eurm10\font\euVII=eurm7\font\euV=eurm5
                                                                \textfont13=\euX
                                                        
        \scriptfont13=\euVII\scriptscriptfont13=\euV
                                                                \def\eu{\fam13\euX}
                                \comm(eu)Passe \`a la famille \<<euler\>> (fontes eurm); {\eu
                                                                Ceci est du eu}!!!

\mathchardef\eup="0D3A
                                \comm(eup)Ecrit $\eup$ (point en fonte \<<euler\>>); mode
                                                                math\'ematique!!!

\mathchardef\euv="0D3B
                                \comm(euv)Ecrit $\euv$ (virgule en fonte \<<euler\>>); mode
                                                                math\'ematique!!!

\let\ev=\emptyset
                                \comm(ev)Ecrit $\ev$; mode math\'ematique!!!

\let\ex=\exists
                                \comm(ex)Ecrit $\ex$; mode math\'ematique!!!

                                \comm(Ex)Provoque un grand saut, puis \'ecrit {\bf Example.-
}!!!

\mathchardef\f="7127
                                \comm(f)Ecrit $\f$; mode math\'ematique!!!

\mathchardef\F="7108
                                \comm(F)Ecrit $\F$; mode math\'ematique!!!
                 
\let\fl=\longrightarrow
                                \comm(fl)Ecrit $\fl$; mode math\'ematique!!!

                                \comm(Fin)Ecrit $\carre$ pr\'ec\'ed\'e d'un blanc inseccable!!!

\let\flc=\rightarrow
                                \comm(flc)Ecrit $\flc$ (fl\`eche courte); mode math\'ematique!!!

{\catcode`@=11
                                                                \catcode`\;=\active%
                                                                \catcode`\:=\active%
                                                                \catcode`\!=\active%
                                                                \catcode`\?=\active%
                                \gdef\francais{
                                                                \catcode`\;=\active%
                                                        
        \def;{\relax\ifhmode\ifdim\lastskip>\z@\unskip\fi
                                                                \kern.2em\fi\string;}
                                                                \catcode`\:=\active%
                                                        
        \def:{\relax\ifhmode\ifdim\lastskip>\z@\unskip\fi
                                                                \kern.2em\fi\string:}
                                                                \catcode`\!=\active%
                                                        
        \def!{\relax\ifhmode\ifdim\lastskip>\z@\unskip\fi
                                                                \kern.2em\fi\string!}
                                                                \catcode`\?=\active%
                                                        
        \def?{\relax\ifhmode\ifdim\lastskip>\z@\unskip\fi
                                                                \kern.2em\fi\string?}
                                                                \frenchspacing\catcode`\@=12}}
                                \comm(francais)Adapte \`a la typographie francaise;
                                                                redefinit les caracteres ; : ? !  aux 
normes francaises, c'est \`a
                                                                dire avec un blanc inseccable devant; 
cela ne change rien de taper
                                                                ou de ne pas taper les blancs avant 
les caracteres ; : ! ? (par contre
                                                                cela change d'en taper apres ou 
non)!!!

\mathchardef\g="710D
                                \comm(g)Ecrit $\g$; mode math\'ematique!!!

\mathchardef\G="7100
                                \comm(G)Ecrit $\G$; mode math\'ematique!!!
                  
%famille 12: go
                                                        
        \font\eufmX=eufm10\font\eufmVII=eufm7\font\eufmV=eufm5
                                                        
        \textfont12=\eufmX\scriptfont12=\eufmVII
                                                                \scriptscriptfont12=\eufmV
                                                                \def\go{\fam12\eufmX}
                                \comm(go)Passe dans la famille \<<gothique\>> (fontes eufm); 
                                                                {\go Ceci est du go}!!!

\mathchardef\gpref="3A40
                                \comm(gpref)Ecrit $\gpref$ (le g signifie \<<grand\>>) ; mode
                                                                math\'ematique!!!
                                                        
\let\gprefe=\sqsubseteq
                                \comm(gprefe)Ecrit $\gprefe$ (le g signifie \<<grand\>>) ; mode 
                                                                math\'ematique!!!

\mathchardef\gsuff="3A41
                                \comm(gsuff)Ecrit $\gsuff$ (le g signifie \<<grand\>>) ; mode
                                                                math\'ematique!!!
                                                        
\let\gsuffe=\sqsupseteq
                                \comm(gsuffe)Ecrit $\gsuffe$ (le g signifie \<<grand\>>) ; mode 
                                                                math\'ematique!!!

\mathchardef\h="7112
                                \comm(h)Ecrit $\h$; mode math\'ematique!!!

\mathchardef\H="7102
                                \comm(H)Ecrit $\H$; mode math\'ematique!!!
                  
\font\hel=cmb10 at 10 pt
               \font\helVII=cmb10 at 7 pt
               \font\helXII=cmb10 at 12 pt
               \font\helXIV=cmb10 at 14 pt
                                \comm(hel)Passe dans la fonte \<<cmb10\>>; 
                                                                {\hel Ceci est du hel};
                                                                suivie \'eventuellemnt d'une taille en 
romain majuscule;
                                                                disponibles: {\tt \\helVII, \\helXII, 
\\helXIV}!!!

\font\helbf=cmb10 at 10 pt
               \font\helbfXII=cmb10 at 12 pt
               \font\helbfXIV=cmb10 at 14 pt
               \font\helbfXVIII=cmb10 at 18 pt
               \font\helbfXXIV=cmb10 at 24 pt
               \font\helbfXXXVI=cmb10 at 36 pt
                                \comm(helbf)Passe \`a la fonte \<<cmb10 gras\>>;
                                                         {\helbf ceci est du helbf};
                                                                suivie \'eventuellemnt d'une taille en 
chiffre romain majuscule; 
                                                                disponibles:
                                                                {\tt\\helbfXII, \\helbfXIV, 
\\helbfXVIII,\\helbfXXIV,
                                                                \\helbfXXXVI}!!!

\font\helbi=cmb10 at 10 pt
               \font\helbiXII=cmb10 at 12 pt
               \font\helbiXIV=cmb10 at 14 pt
               \font\helbiXVIII=cmb10 at 18 pt
               \font\helbiXXIV=cmb10 at 24 pt
               \font\helbiXXXVI=cmb10 at 36 pt
                                \comm(helbi)Passe \`a la fonte \<<cmb10 gras italique\>>; 
                                                                {\helbi Ceci est du helbi}; 
                                                                suivie \'eventuellemnt d'une taille en 
chiffres romains majuscules;
                                                                disponibles:
                                                                {\tt\\helbiXII, \\helbiXIV, 
\\helbiXVIII, \\helbiXXIV, 
                                                                \\helbiXXXVI}!!!

\def\hfl#1#2{\smash{\mathop{\hbox to 12truemm{\rightarrowfill}}
                                                        
        \limits^{\scriptstyle#1}_{\scriptstyle#2}}}
                                \comm(hfl{\#1}{\#2})Trace une fl\`eche horizontale de 12 mm 
avec
                                                                {\tt\#1} dessus et {\tt\#2} dessous!!!

                                \comm(hhat)Ecrit un grand chapeau!!!

\def\ie{\hbox{\it i.e. }}
                                \comm(ie)Ecrit \ie avec un blanc apr\`es!!!

\def\iff{if and only if }
                                \comm(iff)Ecrit \iff avec un blanc apr\`es!!!

\def\Im{{\be Im}}
                                \comm(Im)Ecrit \Im!!!

\let\imp=\Longrightarrow
                                \comm(imp)Ecrit $\imp$; mode math\'ematique; autre nom:
                                                                {\tt\\impl}!!!

\let\impl=\Longrightarrow
                                \comm(impl)Voir {\tt\\imp}!!!

\let\inc=\subset
                                \comm(inc)Ecrit $\inc$; mode math\'ematique!!!

\let\ince=\subseteq
                                \comm(ince)Ecrit $\ince$; mode math\'ematique!!!

\def\inserer(#1)(#2)(#3)(#4)(#5){
                                                                \dimen1=#2\divide\dimen1 
by1000\multiply\dimen1 by#4
                                                                \dimen2=#3\divide\dimen2 
by1000\multiply\dimen2 by#4
                                                        
        \midinsert\vglue\dimen2$$\vbox{\hsize=\dimen1
                                                                \ni\special{picture #1 scaled #4}}$$
                                                                \centerline{#5}\endinsert}
                                \comm(inserer)Ins\`ere au mieux un dessin
                                                                la syntaxe est {\tt\\inserer(nom de la 
figure)(largeur)(hauteur)
                                                                (echelle)(commentaire)}
                                                                (les parametres largeur et hauteur
                                                                sont \`a lire dans le fichier Pictures; 
il faut taper l'unite)!!!

\let\inter=\cap
                                \comm(inter)Ecrit $\inter$; mode math\'ematique!!!

\let\Inter=\bigcap
                                \comm(Inter)Ecrit $\Inter$; mode math\'ematique!!!

%famille 4: it
                                                        
        \font\cmtiX=cmti10\font\cmtiVII=cmti7 
                                                                \font\cmtiV=cmti10 at 5pt 
                                                        
        \textfont4=\cmtiX\scriptfont4=\cmtiVII\scriptscriptfont4=\cmtiV
                                                                \def\it{\fam4\cmtiX}

\mathchardef\k="7114
                                \comm(k)Ecrit $\k$; mode math\'ematique!!!

\mathchardef\l="7115
                                \comm(l)Ecrit $\l$; mode math\'ematique!!!

\mathchardef\L="7103
                                \comm(L)Ecrit $\L$; mode math\'ematique!!!
                  
\def\Lem#1{\bgni{\bf Lemma#1.-}}
                                \comm(Lem)La commande {\tt\\Lem\#} provoque un saut, puis 
\'ecrit
                                                                sans indentation {\bf Lemma\#.-}!!!

                                \comm(lettre)Ecrit le haut de page pour une lettre en anglais; il faut
                                                                le sceau dans le fichier pictures!!!

                                \comm(lettref)Ecrit le haut de page pour une lettre en francais; il 
faut
                                                                le sceau dans le fichier pictures!!!

                                \comm(lettrem)Ecrit le haut de page pour une lettre en francais!!!

\def\lpol{\leavevmode\setbox0=\hbox{l}\hbox to \wd0{\hss\char'40l}}
                                \comm(lpol)Ecrit \lpol\ (l polonais)!!!

\def\Lpol{\leavevmode\setbox0=\hbox{L}\hbox to \wd0{\hss\char'40L}}
                                \comm(Lpol)Ecrit \Lpol\ (L polonais)!!!

\mathchardef\m="7116
                                \comm(m)Ecrit $\m$; mode math\'ematique!!!

                                \comm(mg)Provoque un saut moyen et favorise la coupure!!!

                                \comm(mgni)Provoque un saut moyen, favorise la coupure et 
supprime
                                                                l'indentation!!!

%famille 1: mi
                                                                \def\mi{\fam1\teni}
                                \comm(mi)Passe \`a la famille \<<math. italique\>> (fontes cmmi); 
{\mi
                                                                ceci est du mi}; autres noms: 
{\tt\\mit} (mode math\'ematique),
                                                                {\tt\\oldstyle} (tout mode)!!!

                                \comm(mn)Provoque un saut moyen et d\'efavorise la coupure!!!

                                \comm(mnni)Provoque un saut moyen, d\'efavorise la coupure et
                                                                supprime l'indentation!!!

\def\move#1#2#3
                                                                {\vbox to0pt{\kern-
#3mm\smash{\hbox{\kern#2mm{#1}}}\vss}
                                                                \nointerlineskip}
                                \comm(move\#1\#2\#3)Ecrit \#1 \`a une abscisse de \#2 mm et \`a
                                                                une ordonn\'ee de \#3 mm par 
rapport au point courant (qui n'est
                                                                pas modifi\'e); mode vertical. 
Utilisable pour figures. Fabriquer
                                                                une \\{\tt vbox} suivant la syntaxe 
suivante \par
                                                                \hskip3cm{\tt \\vbox to ... 
mm\{\\hsize=... mm\\vfill\par
                                                         \hskip3cm\\special \{picture ... scaled ... \} 
\par 
                                                         \hskip3cm\\move \{\ ... \}\{\ ... \}\{\ ... \} 
\par 
                                                                \hskip3cm\\hfill\} }\par
                                                                Le {\tt\\vfill} est pour tasser la boite 
vers le bas, le {\tt\\hfill}
                                                                est pour tasser vers la gauche et 
assurer la largeur!!!

\mathchardef\n="7117
                                \comm(n)Ecrit $\n$; mode math\'ematique!!!

\def\Nat{\hbox{\rm I\kern-.9truemm\hbox{N}}}
                                \comm(Nat)Ecrit \Nat!!!

\let\ni=\noindent
                                \comm(ni)Supprime l'indentation!!!

\mathchardef\NN="0B4E
                                \comm(NN)Ecrit $\NN$: mode math\'ematique!!!

\def\non{\mathord{\be non}}
                                \comm(non)Ecrit $\non$; mode math\'ematique!!!

                                \comm(notapp)Ecrit $\notin$ sans blanc avant ni apr\`es; mode
                                                                        math\'ematique!!!

\mathchardef\o="7121
                                \comm(o)Ecrit $\o$; mode math\'ematique!!!

\mathchardef\O="710A
                                \comm(O)Ecrit $\O$; mode math\'ematique!!!
                  
\let\oo=\infty
                                \comm(oo)Ecrit $\oo$; mode math\'ematique!!!

\def\ou{\hbox{ ou }}
                                \comm(ou)Ecrit $\ou$ en tout mode (avec un blanc avant et 
apr\`es)!!!

\mathchardef\p="7119
                                \comm(p)Ecrit $\p$; mode math\'ematique!!!

\mathchardef\P="7105
                                \comm(P)Ecrit $\P$; mode math\'ematique!!!

\def\pmb#1{\setbox0=\hbox{#1}\kern-.15pt\copy0\kern-\wd0
                 \kern-.3pt\copy0\kern-\wd0\kern-.15pt\raise.1pt\box0}
                                \comm(pmb)Met en gras son param\`etre: {\tt\\pmb\#} \'ecrit \# en
                                                                faux gras (ce qui donne 
\pmb{\#})!!!

\def\pp{\hbox{$\ldots$}}
                                \comm(pp)Ecrit \pp!!!

\def\Ppar{\mathhexbox27B}
                                \comm(Ppar)Ecrit \Ppar!!!

\def\ppp{\hbox{$, \ldots, $}}
                                \comm(ppp)Ecrit \ppp!!!

\let\pr=\vdash
                                \comm(pr)Ecrit $\pr$; mode math\'ematique!!!

\def\Pr{\bgni{\it Proof. }}
                                \comm(Pr)Fait un saut et \'ecrit \<<{\it Proof.}\>>!!!

\def\pref{\mathrel{\scriptstyle\gpref}}
                                \comm(pref)Ecrit la relation $\pref$; mode math\'ematique!!!
                                                        
\def\prefe{\mathrel{\scriptstyle\gprefe}}
                                \comm(prefe)Ecrit la relation $\prefe$; mode math\'ematique!!!
                                                        
\def\Prop#1{\bgni{\bf Proposition#1.-}}
                                \comm(Prop)La commande {\tt\\Prop\#} provoque un saut, puis 
\'ecrit
                                                                sans indentation {\bf Proposition\#.-
}!!!

\mathchardef\Pw="0C50
                                \comm(Pw)Ecrit $\Pw$; mode math\'ematique!!!

\mathchardef\q="711F
                                \comm(q)Ecrit $\q$; mode math\'ematique!!!

\let\qq=\forall
                                \comm(qq)Ecrit $\qq$; mode math\'ematique!!!

\mathchardef\r="711A
                                \comm(r)Ecrit $\r$; mode math\'ematique!!!

\mathchardef\res="0A16
                                \comm(res)Ecrit $\res$; mode math\'ematique!!!
                                                        
\def\resp{\hbox{\it resp. }}
                                \comm(resp)Ecrit \resp (avec un blanc apr\`es)!!!

\def\Ree{\hbox{\rm I\kern-.9truemm\hbox{R}}}
                                \comm(Ree)Ecrit \Ree!!!

\mathchardef\RR="0B52
                                \comm(RR)Ecrit $\RR$; mode math\'ematique!!!

\def\rres{\hbox to 4.5pt{$\res$\kern-1pt\hbox{$\res$}\hss}}
                                \comm(rres)Ecrit $\rres$; mode math\'ematique!!!

\mathchardef\s="711B
                                \comm(s)Ecrit $\s$; mode math\'ematique!!!

\mathchardef\S="7106
                                \comm(S)Ecrit $\S$; mode math\'ematique!!!

%famille10: sa;
                                                        
        \font\msamX=msam10\font\msamVII=msam7 
                                                                \font\msamV=msam5 
\textfont10=\msamX
                                                        
        \scriptfont10=\msamVII\scriptscriptfont10=\msamV
                                                                \def\sa{\fam10\msamX}
                                \comm(sa)Passe \`a la famille \<<symboles a\>> (fontes msam); 
{\sa
                                                                Ceci est du sa}!!!

%famille11: sb;
                                                        
        \font\msbmX=msbm10\font\msbmVII=msbm7 
                                                                \font\msbmV=msbm5 
\textfont11=\msbmX
                                                        
        \scriptfont11=\msbmVII\scriptscriptfont11=\msbmV
                                                                \def\sb{\fam11\msbmX}
                                                                \let\db=\sb
                                \comm(sb)Passe \`a la famille \<<symboles b\>> (fontes msbm); 
{\sb
                                                                CECI EST DU SB}; autre nom: 
{\tt\\db}!!!

\let\sat=\models
                                \comm(sat)Ecrit $\sat$; mode math\'ematique!!!

                                \comm(sg)Provoque un saut petit et favorise la coupure!!!

\def\sgni{\smallbreak\noindent}
                                \comm(sgni)Provoque un saut petit, favorise la coupure et 
supprime
                                                                l'indentation!!!

%famille 5: sl
                                                        
        \font\cmslX=cmsl10\font\cmslVII=cmsl10 at 7 pt 
                                                                \font\cmslV=cmsl10 at 5pt 
                                                        
        \textfont5=\cmslX\scriptfont5=\cmslVII\scriptscriptfont5=\cmslV
                                                                \def\sl{\fam5\cmslX}

\def\sland{\hbox{ \sl and }}
                                \comm(sland)Ecrit $\sland$ en tout mode!!!

\font\smcap=cmcsc10                                                             
                                \comm(smcap)Passe \`a la fonte \<<petites capitales\>>\ (cmcsc);
                                                                {\smcap Ceci est du smcap}!!!

                                \comm(sn)Provoque un saut petit et d\'efavorise la coupure!!!

                                \comm(snni)Provoque un saut petit, d\'efavorise la coupure et
                                                                supprime l'indentation!!!

\def\Spar{\mathhexbox278}
                                \comm(Spar)Ecrit \Spar!!!

\mathchardef\SS="0B53
                                \comm(SS)Ecrit $\SS$: mode math\'ematique!!!

\def\ssi{si et seulement si }
                                \comm(ssi)Ecrit \ssi (avec un blanc apr\`es)!!!

\def\suff{\mathrel{\scriptstyle\gsuff}}
                                \comm(suff)Ecrit la relation $\suff$; mode math\'ematique!!!
                                                        
\def\suffe{\mathrel{\scriptstyle\gsuffe}}
                                \comm(suffe)Ecrit la relation $\suffe$; mode math\'ematique!!!
                                                        
\def\Supp{\mathord{\be Supp}}
                                \comm(Supp)Ecrit $\Supp$; mode math\'ematique!!!

%famille 2: sy
                                                                \def\sy{\fam2\tensy}
                                                                \let\ca=\sy
                                \comm(sy)Passe \`a la famille \<<symbole\>> (fontes cmsy); {\sy
                                                                CECI EST DU SY}; autres noms: 
{\tt\\cal} (mode math\'ematique),
                                                                {\tt\\ca} (tout mode)!!!

\mathchardef\t="711C
                                \comm(t)Ecrit $\t$; mode math\'ematique!!!

\def\Thm#1{\bgni{\bf Theorem#1.-}}
                                \comm(Thm)La commande {\tt\\Thm\#} provoque un saut, puis 
\'ecrit
                                                                sans indentation {\bf Theorem\#.-
}!!!

\font\tim=cmb10 at 10pt 
                
               \font\timXIV=cmb10 at 14pt 
               \font\timXVIII=cmb10 at 18pt 
               \font\timVII=cmb10 at 7pt 
                                \comm(tim)Passe dans la fonte \og cmb10\fg; 
                                                                {\tim Ceci est du tim};
                                                                suivie \'eventuellemnt d'une taille en 
chiffres romains majuscules;
                                                                disponibles: {\tt\\timVII, \\timlXII, 
\\timXIV, \\timXVIII}!!!

\font\timbf=cmb10 at 10pt 
               
              \font\timbfXIV=cmb10 at 14pt 
              \font\timbfXVIII=cmb10 at 18pt 
              \font\timbfVII=cmb10 at 7pt
                                \comm(timbf)Passe dans la fonte \og cmb10 gras\fg; 
                                                                {\timbf Ceci est du timbf};
                                                                suivie \'eventuellemnt d'une taille en 
chiffres romains majuscules;
                                                                disponibles: {\tt\\timbfVII, 
\\timbflXII, \\timbfXIV,
                                                                \\timbfXVIII}!!!

%famille 7: tt
                                                        
        \font\cmttX=cmtt10\font\cmttVII=cmtt10 at 7 pt 
                                                                \font\cmttV=cmtt10 at 5pt 
                                                        
        \textfont7=\cmtiX\scriptfont7=\cmttVII\scriptscriptfont7=\cmttV
                                                                \def\tt{\fam7\cmttX}

                                \comm(ttil)Ecrit un grand tilde!!!

\mathchardef\u="711D
                                \comm(u)Ecrit $\u$; mode math\'ematique!!!

\mathchardef\U="7107
                                \comm(U)Ecrit $\U$; mode math\'ematique!!!

\let\un=\cup
                                \comm(un)Ecrit $\un$; mode math\'ematique!!!

\let\Un=\bigcup
                                \comm(Un)Ecrit $\Un$; mode math\'ematique!!!

\def\val{\mathop{\be val}}
                                \comm(val)Ecrit $\val$; mode math\'ematique!!!

\def\var{{\be var}}
                                \comm(var)Ecrit $\var$!!!

                                \comm(vfl{\#1}{\#2})Trace une fl\`eche verticale de 12 mm avec
                                                                {\tt\#1} \`a gauche et {\tt\#2} \`a 
droite!!!

\mathchardef\w="7120
                                \comm(w)Ecrit $\w$; mode math\'ematique!!!

\mathchardef\W="7109
                                \comm(W)Ecrit $\W$; mode math\'ematique!!!

                                \comm(wrt)Ecrit \<<with respect to\>> suivi d'un blanc 
inseccable!!!

\mathchardef\x="7118
                                \comm(x)Ecrit $\x$; mode math\'ematique!!!

\mathchardef\X="7104
                                \comm(X)Ecrit $\X$; mode math\'ematique!!!

\mathchardef\y="7111
                                \comm(y)Ecrit $\y$; mode math\'ematique!!!

\mathchardef\z="7110
                                \comm(z)Ecrit $\z$; mode math\'ematique!!!

\def\ZZ{{\db Z}}
                                \comm(ZZ)Ecrit $\ZZ$!!!

                                \comm( )\vfill\eject\centerline{\bfXII Anciennes commandes}!!!

\let\bfe=\be
                                \comm(bfe)Ancien nom de {\tt\\be}!!!

                                \comm(bfi)Ancien nom de {\tt\\bi}!!!
  
\let\ca=\sy
                                \comm(ca)Ancien nom de la famille {\tt\\sy}!!!

\font\cmmibX=cmmib10
                                \comm(cmmibX)Passe \`a la fonte cmmib10; {\cmmibX Ceci est 
du
                                                                cmmib}!!!

\let\db=\sb
                                \comm(db)Ancien nom de la famille {\tt\\sb}!!!

\def\fg{~{\cyr\char'076}}
                                \comm(fg)Ancien nom de {\tt\\>>}!!!

                                \comm(msymX)Passe \`a la fonte \<<msbm10\>>; obsol\`ete: 
nouveau
                                                                nom de famille: {\tt\\db}!!!

\def\og{{\cyr\char'074}~}
                                \comm(og)Ancien nom de {\tt\\<<}!!!

                                \comm(rd)Ancien nom de la famille {\tt\\mi}!!!

\def\picture #1 by #2 (#3){
 \vbox to #2{\hrule width #1 height 0pt depth 0pt
 \vfill\special{picture #3}}}
\def\dessin(#1)(#2)(#3)(#4){{
 \dimen0=#2 \dimen1=#3
 \divide\dimen0 by 1000 \multiply\dimen0 by #4
 \divide\dimen1 by 1000 \multiply\dimen1 by #4
 \picture \dimen0 by \dimen1 (#1 scaled #4)}}

\newcount\numlem
\newcount\numchap
\hsize=12cm
\vsize=20truecm
\voffset=1cm
\hoffset=5mm
\newtoks\gauche \gauche={Serge Burckel}
\newtoks\droite 
\def\makeheadline{\vbox to 0pt{\vskip -40pt\line{\vbox to 8.5pt{}\the
\headline}\vss\nointerlineskip}}
\headline={{\voffset 2cm\sevenbf\the\gauche\hfill\the\droite}}

\def\sk{\smallskip}

\overfullrule=0mm 
\parindent=0mm

%\catcode`\-=\active
%\def-{\ifmmode{\!\hskip 0.7pt\char 45\!\hskip
%0.7pt}\else{\char 45}\fi}
 
%\catcode`\+=\active
%\def+{\ifmmode{\!\hskip 0.7pt\char 43\!\hskip
%0.7pt}\else{\char 43}\fi}

\def\ca{{\bi a}}

\def\l{n}

\def\psi{{\bf cl}}
\def\x{X}
\def\fl{{\longrightarrow}}

\def\H#1{{\hskip #1truemm}}
\def\BX#1{\boxit{8pt}{\vbox{#1}}}

\def\s#1 #2{{#1}_{(#2)}}
\def\O{{\sy T}}
\def\C#1 #2{{\(#1\)_{#2}}}

\def\D{\L}

\def\ttchap#1{
\numlem=0
\vfill\eject\gauche={\the\numchap.
#1}\droite={\hfill}{\bfXVIII \BX{CHAPITRE 
\the\numchap. \bn \hbox{#1}} }\bn\vskip 2cm}

\def\ttparag
#1{{\bfXII #1}\bn}

 \def\Pr{\bg\bf PREUVE.
\rm} \def\Rp{\hbox{\carre\rm}\bn}

\def\Def#1#2{\bg{{\bf Definition. #1}#2}\rm\bg}

\def\Lem#1#2{\advance\numlem by
1\bg{{\bf LEMME \the\numlem. #1}\sl#2}\rm\bg}
\def\Lemprop#1#2{\advance\numlem by 1\bg{{\bf
PROPOSITION \the\numlem. #1}\sl#2}\rm\bg}
\def\Lempropr#1#2{\advance\numlem by 1\bg{{\bf
PROPRIETE \the\numlem. #1}\sl#2}\rm\bg}
\def\Thm#1#2{\advance\numlem by 1\bg{{\bf
THEOREME \the\numlem. #1}\sl#2}\rm\bg}

\font\sc=cmcsc10

\def\ttchap#1{
\numlem=0
\vfill\eject\gauche={\the\numchap.
#1}\droite={\hfill}{\bfXVIII \BX{CHAPTER  \the\numchap.\bn
\hbox{#1}} }\bn\vskip 2cm}

\def\Pr{\bg\bf Proof.
\rm}

\def\Def#1#2{\bg{{\bf Definition. #1}#2}\rm\bg}

\def\Lem#1#2{\advance\numlem by
1\bg{{\bf Lemma \the\numlem. #1}\sl#2}\rm\bg}
\def\Lemprop#1#2{\advance\numlem by 1\bg{{\bf
Proposition \the\numlem. #1}\sl#2}\rm\bg}
\def\Prop#1#2{\bg{{\bf
Proposition  #1}\sl#2}\rm\bg}
\def\Lempropr#1#2{\advance\numlem by 1\bg{{\bf
Property \the\numlem. #1}\sl#2}\rm\bg}
\def\Thm#1#2{\advance\numlem by 1\bg{{\bf
Theorem \the\numlem. #1}\sl#2}\rm\bg}

\catcode`\Ž=\active \def Ž{\'e}
\catcode`\ˆ=\active \def ˆ{\`a}
\catcode`\=\active \def {\`u}
\catcode`\=\active \def {\`e}
\catcode`\=\active \def {\^e}
\catcode`\™=\active \def ™{\^o}
\catcode`\"=\active \def "{\^i}
\catcode`\=\active \def {\c c}

\def\fileversion{v1.0}
\def\filedate{3 Dec 91}
\immediate\write16{Document style option `epsfig', \fileversion\space
<\filedate> (edited by SPQR)}
\chardef\atcode=\catcode`\@
\catcode`\@=11\relax
\ifx\typeout\undefined%
        \newwrite\@unused
        \def\typeout#1{{\let\protect\string\immediate\write\@unused{#1}}}
\fi
\ifx\undefined\epsfig%
\else
        \typeout{EPSFIG --- already loaded} 
\fi
%
%
% we try to put a box round space of missing figures. plain TeX
% doesn't have an easy command, so just ignore it
\ifx\undefined\fbox\def\fbox#1{#1}\fi
%%%
%%% we need Rokicki's EPSF macros anyway, unless they are already loaded
%
\ifx\undefined\epsfbox\input epsf\fi
%
%% SPQR 12.91 handling of errors using standard LaTeX error 
%% mechanism. In case we are plain TeX we first define the
%% error routines...
\ifx\undefined\@latexerr
        \newlinechar`\^^J
        \def\@spaces{\space\space\space\space}
        \def\@latexerr#1#2{%
        \edef\@tempc{#2}\expandafter\errhelp\expandafter{\@tempc}%
        \typeout{Error. \space see a manual for explanation.^^J
         \space\@spaces\@spaces\@spaces Type \space H <return> \space for
         immediate help.}\errmessage{#1}}
\fi
%------------------------
%% a couple of LaTeX error messages
\def\@whattodo{You tried to include a PostScript figure which 
cannot be found^^JIf you press return to carry on anyway,^^J
The failed name will be printed in place of the figure.^^J
or type X to quit}
\def\@whattodobb{You tried to include a PostScript figure which 
has no^^Jbounding box, and you supplied none.^^J
If you press return to carry on anyway,^^J
The failed name will be printed in place of the figure.^^J
or type X to quit}
%------------------------
%
\newwrite\@unused
%------------------------------------------------------------------------
%------------------------------------------------------------------------
%%% @psdo control structure -- similar to Latex @for.
%%% I redefined these with different names so that psfig can
%%% be used with TeX as well as LaTeX, and so that it will not 
%%% be vunerable to future changes in LaTeX's internal
%%% control structure,
%
\def\@nnil{\@nil}
\def\@empty{}
\def\@psdonoop#1\@@#2#3{}
\def\@psdo#1:=#2\do#3{\edef\@psdotmp{#2}\ifx\@psdotmp\@empty \else
    \expandafter\@psdoloop#2,\@nil,\@nil\@@#1{#3}\fi}
\def\@psdoloop#1,#2,#3\@@#4#5{\def#4{#1}\ifx #4\@nnil \else
       #5\def#4{#2}\ifx #4\@nnil \else#5\@ipsdoloop #3\@@#4{#5}\fi\fi}
\def\@ipsdoloop#1,#2\@@#3#4{\def#3{#1}\ifx #3\@nnil 
       \let\@nextwhile=\@psdonoop \else
      #4\relax\let\@nextwhile=\@ipsdoloop\fi\@nextwhile#2\@@#3{#4}}
\def\@tpsdo#1:=#2\do#3{\xdef\@psdotmp{#2}\ifx\@psdotmp\@empty \else
    \@tpsdoloop#2\@nil\@nil\@@#1{#3}\fi}
\def\@tpsdoloop#1#2\@@#3#4{\def#3{#1}\ifx #3\@nnil 
       \let\@nextwhile=\@psdonoop \else
      #4\relax\let\@nextwhile=\@tpsdoloop\fi\@nextwhile#2\@@#3{#4}}
%%% 
%
%%%%%%%%%%%%%%%%%%%%%%%%%%%%%%%%%%%%%%%%%%%%%%%%%%%%%%%%%%%%%%%%%%%
%%% file reading stuff from epsf.tex
%%%   EPSF.TEX macro file:
%%%   Written by Tomas Rokicki of Radical Eye Software, 29 Mar 1989.
%%%   Revised by Don Knuth, 3 Jan 1990.
%%%   Revised by Tomas Rokicki to accept bounding boxes with no
%%%      space after the colon, 18 Jul 1990.
%%%   Portions modified/removed for use in PSFIG package by
%%%      J. Daniel Smith, 9 October 1990.
%%%   Just the bit which knows about (atend) as a BoundingBox
%
%%%    hacked back a bit by SPQR 12/91
%
\long\def\epsfaux#1#2:#3\\{\ifx#1\epsfpercent
   \def\testit{#2}\ifx\testit\epsfbblit
        \@atendfalse
        \epsf@atend #3 . \\%
        \if@atend       
           \if@verbose
                \typeout{epsfig: found `(atend)'; continuing search}
           \fi
        \else
                \epsfgrab #3 . . . \\%
                \global\no@bbfalse
        \fi
   \fi\fi}%
%
%%% Determine if the stuff following the %%BoundingBox is `(atend)'
%%% J. Daniel Smith.  Copied from \epsf@grab above.
%
\def\epsf@atendlit{(atend)} 
\def\epsf@atend #1 #2 #3\\{%
   \def\epsf@tmp{#1}\ifx\epsf@tmp\empty
      \epsf@atend #2 #3 .\\\else
   \ifx\epsf@tmp\epsf@atendlit\@atendtrue\fi\fi}

%%% End of file reading stuff from epsf.tex
%%%%%%%%%%%%%%%%%%%%%%%%%%%%%%%%%%%%%%%%%%%%%%%%%%%%%%%%%%%%%%%%%%%

%%%%%%%%%%%%%%%%%%%%%%%%%%%%%%%%%%%%%%%%%%%%%%%%%%%%%%%%%%%%%%%%%%%
%%% trigonometry stuff from "trig.tex"
\chardef\letter = 11
\chardef\other = 12

\newif \ifdebug %%% turn me on to see TeX hard at work ...
\newif\ifc@mpute %%% don't need to compute some values
\newif\if@atend
\c@mputetrue % but assume that we do

\let\then = \relax
\def\r@dian{pt }
\let\r@dians = \r@dian
\let\dimensionless@nit = \r@dian
\let\dimensionless@nits = \dimensionless@nit
\def\internal@nit{sp }
\let\internal@nits = \internal@nit
\newif\ifstillc@nverging
\def \Mess@ge #1{\ifdebug \then \message {#1} \fi}

{ %%% Things that need abnormal catcodes %%%
        \catcode `\@ = \letter
        \gdef \nodimen {\expandafter \n@dimen \the \dimen}
        \gdef \term #1 #2 #3%
               {\edef \t@ {\the #1}%%% freeze parameter 1 (count, by value)
                \edef \t@@ {\expandafter \n@dimen \the #2\r@dian}%
                                   %%% freeze parameter 2 (dimen, by value)
                \t@rm {\t@} {\t@@} {#3}%
               }
        \gdef \t@rm #1 #2 #3%
               {{%
                \count 0 = 0
                \dimen 0 = 1 \dimensionless@nit
                \dimen 2 = #2\relax
                \Mess@ge {Calculating term #1 of \nodimen 2}%
                \loop
                \ifnum  \count 0 < #1
                \then   \advance \count 0 by 1
                        \Mess@ge {Iteration \the \count 0 \space}%
                        \Multiply \dimen 0 by {\dimen 2}%
                        \Mess@ge {After multiplication, term = \nodimen 0}%
                        \Divide \dimen 0 by {\count 0}%
                        \Mess@ge {After division, term = \nodimen 0}%
                \repeat
                \Mess@ge {Final value for term #1 of 
                                \nodimen 2 \space is \nodimen 0}%
                \xdef \Term {#3 = \nodimen 0 \r@dians}%
                \aftergroup \Term
               }}
        \catcode `\p = \other
        \catcode `\t = \other
        \gdef \n@dimen #1pt{#1} %%% throw away the ``pt''
}

\def \Divide #1by #2{\divide #1 by #2} %%% just a synonym

\def \Multiply #1by #2%%% allows division of a dimen by a dimen
       {{%%% should really freeze parameter 2 (dimen, passed by value)
        \count 0 = #1\relax
        \count 2 = #2\relax
        \count 4 = 65536
        \Mess@ge {Before scaling, count 0 = \the \count 0 \space and
                        count 2 = \the \count 2}%
        \ifnum  \count 0 > 32767 %%% do our best to avoid overflow
        \then   \divide \count 0 by 4
                \divide \count 4 by 4
        \else   \ifnum  \count 0 < -32767
                \then   \divide \count 0 by 4
                        \divide \count 4 by 4
                \else
                \fi
        \fi
        \ifnum  \count 2 > 32767 %%% while retaining reasonable accuracy
        \then   \divide \count 2 by 4
                \divide \count 4 by 4
        \else   \ifnum  \count 2 < -32767
                \then   \divide \count 2 by 4
                        \divide \count 4 by 4
                \else
                \fi
        \fi
        \multiply \count 0 by \count 2
        \divide \count 0 by \count 4
        \xdef \product {#1 = \the \count 0 \internal@nits}%
        \aftergroup \product
       }}

\def\r@duce{\ifdim\dimen0 > 90\r@dian \then   % sin(x) = sin(180-x)
                \multiply\dimen0 by -1
                \advance\dimen0 by 180\r@dian
                \r@duce
            \else \ifdim\dimen0 < -90\r@dian \then  % sin(x) = sin(360+x)
                \advance\dimen0 by 360\r@dian
                \r@duce
                \fi
            \fi}

\def\Sine#1%
       {{%
        \dimen 0 = #1 \r@dian
        \r@duce
        \ifdim\dimen0 = -90\r@dian \then
           \dimen4 = -1\r@dian
           \c@mputefalse
        \fi
        \ifdim\dimen0 = 90\r@dian \then
           \dimen4 = 1\r@dian
           \c@mputefalse
        \fi
        \ifdim\dimen0 = 0\r@dian \then
           \dimen4 = 0\r@dian
           \c@mputefalse
        \fi
        \ifc@mpute \then
                % convert degrees to radians
                \divide\dimen0 by 180
                \dimen0=3.141592654\dimen0
                \dimen 2 = 3.1415926535897963\r@dian %%% a well-known constant
                \divide\dimen 2 by 2 %%% we only deal with -pi/2 : pi/2
                \Mess@ge {Sin: calculating Sin of \nodimen 0}%
                \count 0 = 1 %%% see power-series expansion for sine
                \dimen 2 = 1 \r@dian %%% ditto
                \dimen 4 = 0 \r@dian %%% ditto
                \loop
                        \ifnum  \dimen 2 = 0 %%% then we've done
                        \then   \stillc@nvergingfalse 
                        \else   \stillc@nvergingtrue
                        \fi
                        \ifstillc@nverging %%% then calculate next term
                        \then   \term {\count 0} {\dimen 0} {\dimen 2}%
                                \advance \count 0 by 2
                                \count 2 = \count 0
                                \divide \count 2 by 2
                                \ifodd  \count 2 %%% signs alternate
                                \then   \advance \dimen 4 by \dimen 2
                                \else   \advance \dimen 4 by -\dimen 2
                                \fi
                \repeat
        \fi             
                        \xdef \sine {\nodimen 4}%
                        %\typeout {Sin: sine of #1 \space is \sine \space}%
       }}

%%% Now the Cosine can be calculated easily by calling \Sine:
%%%  cos(x) = sin(90-x)
\def\Cosine#1{\ifx\sine\UnDefined\edef\Savesine{\relax}\else
                             \edef\Savesine{\sine}\fi
        {\dimen0=#1\r@dian\multiply\dimen0 by -1
         \advance\dimen0 by 90\r@dian
         \Sine{\nodimen 0}
         \xdef\cosine{\sine}
         %\typeout {Cosine: cos of \space \nodimen 0 \space is \cosine \space}%
         \xdef\sine{\Savesine}}}              
%%% end of trig stuff
%%%%%%%%%%%%%%%%%%%%%%%%%%%%%%%%%%%%%%%%%%%%%%%%%%%%%%%%%%%%%%%%%%%%
%
\def\psdraft{\def\@psdraft{0}}
\def\psfull{\def\@psdraft{100}}
\psfull
\newif\if@draftbox
\def\psnodraftbox{\@draftboxfalse}
\@draftboxtrue
\newif\if@noisy
\def\pssilent{\@noisyfalse}
\def\psnoisy{\@noisytrue}
\@noisyfalse

%%% These are for the option list.
%%% A specification of the form a = b maps to calling \@p@@sa{b}
\newif\if@bbllx
\newif\if@bblly
\newif\if@bburx
\newif\if@bbury
\newif\if@height
\newif\if@width
\newif\if@rheight
\newif\if@rwidth
\newif\if@angle
\newif\if@clip
\newif\if@verbose
\def\@p@@sclip#1{\@cliptrue}
\newif\if@prologfile
\@prologfiletrue
%
%%% if this is true, the original Darrell macros and specials are used
\newif\ifuse@psfig
\def\@p@@sfile#1{%
\def\@p@sfile{NO FILE: #1}%
\def\@p@sfilefinal{NO FILE: #1}%
        \openin1=#1
        \ifeof1\closein1
                \openin1=#1.bb
                        \ifeof1\closein1
                                \if@bbllx\if@bblly\if@bburx\if@bbury% added 6/91 Rob Russell
                                        \def\@p@sfile{#1}%
                                        \def\@p@sfilefinal{#1}%
                                        \fi\fi\fi
                                \else
                                        \@latexerr{ERROR! PostScript file #1 not found}\@whattodo
                                        \@p@@sbbllx{100bp}
                                        \@p@@sbblly{100bp}
                                        \@p@@sbburx{200bp}
                                        \@p@@sbbury{200bp}
                                        \def\@p@scost{200}
                                \fi
                        \else
                                \closein1%
                                \edef\@p@sfile{#1.bb}%
                                \edef\@p@sfilefinal{"`zcat `texfind #1.Z`"}%
                        \fi
        \else\closein1
                    \edef\@p@sfile{#1}%
                    \edef\@p@sfilefinal{#1}%
        \fi%
}
 % alternative syntax: figure=
\let\@p@@sfigure\@p@@sfile
\def\@p@@sbbllx#1{
                %\typeout{bbllx is #1}
                \use@psfigtrue
                \@bbllxtrue
                \dimen100=#1
                \edef\@p@sbbllx{\number\dimen100}
}
\def\@p@@sbblly#1{
                %\typeout{bblly is #1}
                \use@psfigtrue
                \@bbllytrue
                \dimen100=#1
                \edef\@p@sbblly{\number\dimen100}
}
\def\@p@@sbburx#1{
                %\typeout{bburx is #1}
                \use@psfigtrue
                \@bburxtrue
                \dimen100=#1
                \edef\@p@sbburx{\number\dimen100}
}
\def\@p@@sbbury#1{
                %\typeout{bbury is #1}
                \use@psfigtrue
                \@bburytrue
                \dimen100=#1
                \edef\@p@sbbury{\number\dimen100}
}
\def\@p@@sheight#1{
                \@heighttrue
                \epsfysize=#1
                \dimen100=#1
                \edef\@p@sheight{\number\dimen100}
                %\typeout{Height is \@p@sheight}
}
\def\@p@@swidth#1{
                %\typeout{Width is #1}
                \@widthtrue
                \epsfxsize=#1
                \dimen100=#1
                \edef\@p@swidth{\number\dimen100}
}
\def\@p@@srheight#1{
                %\typeout{Reserved height is #1}
                \@rheighttrue
                \dimen100=#1
                \edef\@p@srheight{\number\dimen100}
}
\def\@p@@srwidth#1{
                %\typeout{Reserved width is #1}
                \@rwidthtrue
                \dimen100=#1
                \edef\@p@srwidth{\number\dimen100}
}
\def\@p@@sangle#1{
                %\typeout{Rotation is #1}
                \use@psfigtrue
                \@angletrue
%               \dimen100=#1
                \edef\@p@sangle{#1} %\number\dimen100}
}
\def\@p@@ssilent#1{ 
                \@verbosefalse
}
\def\@cs@name#1{\csname #1\endcsname}
\def\@setparms#1=#2,{\@cs@name{@p@@s#1}{#2}}
%
%%% initialize the defaults (size the size of the figure)
%
\def\ps@init@parms{
                \@bbllxfalse \@bbllyfalse
                \@bburxfalse \@bburyfalse
                \@heightfalse \@widthfalse
                \@rheightfalse \@rwidthfalse
                \def\@p@sbbllx{}\def\@p@sbblly{}
                \def\@p@sbburx{}\def\@p@sbbury{}
                \def\@p@sheight{}\def\@p@swidth{}
                \def\@p@srheight{}\def\@p@srwidth{}
                \def\@p@sangle{0}
                \def\@p@sfile{}
                \def\@p@scost{10}
                \use@psfigfalse
                \def\@sc{}
                \if@noisy
                        \@verbosetrue
                \else
                        \@verbosefalse
                \fi
                \@clipfalse
}
%
%%% Go through the options setting things up.
%
\def\parse@ps@parms#1{
                \@psdo\@psfiga:=#1\do
                   {\expandafter\@setparms\@psfiga,}}
%
%%% Compute bb height and width
%
\newif\ifno@bb
\def\bb@missing{
        \epsfgetbb{\@p@sfile}
        \ifepsfbbfound\no@bbfalse\else\no@bbtrue\bb@cull\epsfllx\epsflly\epsfurx\epsfury\fi
}       
\def\bb@cull#1#2#3#4{
        \dimen100=#1 bp\edef\@p@sbbllx{\number\dimen100}
        \dimen100=#2 bp\edef\@p@sbblly{\number\dimen100}
        \dimen100=#3 bp\edef\@p@sbburx{\number\dimen100}
        \dimen100=#4 bp\edef\@p@sbbury{\number\dimen100}
        \no@bbfalse
}

\newdimen\p@intvaluex
\newdimen\p@intvaluey
\newdimen\@ffsetvalue
\newdimen\x@ffsetvalue
\newdimen\y@ffsetvalue

%%% Calculate \@ffsetvalue = (#2 - #1) \sin\theta
%%%  The sine of the angle is already stored in \sine.
%%%  If (#2-#1)>0, then the result is zero in the 2nd and 4th quadrants, and
%%%  if (#2-#1)<0, then the result is zero in the 1st and 3rd quadrants.
%%%  Only the x coordinate needs an offset in the 1st and 3rd quadrants,
%%%  and only the y coordinate needs an offset otherwise.

\def\compute@offset#1#2{{\dimen0=#1 sp\dimen1=#2 sp
                        \advance\dimen1 by -\dimen0
                        \dimen1=\sine\dimen1
                        \dimen0=\cosine\dimen1
                        \ifdim\dimen0<0sp \dimen1=0sp \fi
                        \global\@ffsetvalue=\dimen1}}

%%% rotate point (#1,#2) about (0,0).
%%% The sine and cosine of the angle are already stored in \sine and
%%% \cosine.  The result is placed in (\p@intvaluex, \p@intvaluey).
\def\rotate@#1#2{{\dimen0=#1 sp\dimen1=#2 sp
%%%             calculate x' = x \cos\theta - y \sin\theta
                  \global\p@intvaluex=\cosine\dimen0
                  \dimen3=\sine\dimen1
                  \global\advance\p@intvaluex by -\dimen3
%%%             calculate y' = x \sin\theta + y \cos\theta
                  \global\p@intvaluey=\sine\dimen0
                  \dimen3=\cosine\dimen1
                  \global\advance\p@intvaluey by \dimen3
                  }}
%%% rotate point (#1,#2) about the point (#3,#4), finding the x value.
%%% The sine and cosine of the angle are already stored in \sine and
%%% \cosine.  The result is placed in \p@intvaluex
%\def\rotate@x#1#2#3#4{{\dimen0=#1 sp
%                       \dimen1=#2 sp
%                       \dimen2=#3 sp
%                       \dimen4=#4 sp
%                       \advance\dimen0 by -\dimen3
%                       \dimen0=\cosine\dimen0
%                       \advance\dimen4 by -\dimen2
%                       \dimen4=\sine\dimen4
%                       \global\p@intvaluex=\dimen0
%                       \global\advance\p@intvaluex by \dimen4
%                       \global\advance\p@intvaluex by \dimen3
%
%}}
\def\compute@bb{
                \no@bbfalse
                \if@bbllx \else \no@bbtrue \fi
                \if@bblly \else \no@bbtrue \fi
                \if@bburx \else \no@bbtrue \fi
                \if@bbury \else \no@bbtrue \fi
                \ifno@bb \bb@missing \fi
                \ifno@bb
                        \@latexerr{ERROR! cannot locate BB!}\@whattodobb
                        \@p@@sbbllx{100bp}
                        \@p@@sbblly{100bp}
                        \@p@@sbburx{200bp}
                        \@p@@sbbury{200bp}
                        \def\@p@scost{200}
                \fi
                %\typeout{BB: \@p@sbbllx, \@p@sbblly, \@p@sbburx, \@p@sbbury} 
                \if@angle 
                        \Sine{\@p@sangle}\Cosine{\@p@sangle}
                        \compute@offset{\@p@sbblly}{\@p@sbbury}
                        \x@ffsetvalue=\@ffsetvalue
                        % Note that arguments are reversed to
                        %  give a negative interval:
                        \compute@offset{\@p@sbburx}{\@p@sbbllx}
                        \y@ffsetvalue=\@ffsetvalue

                        \rotate@{\@p@sbbllx}{\@p@sbblly}
                        \advance\p@intvaluex by -\x@ffsetvalue
                        \advance\p@intvaluey by -\y@ffsetvalue
                        \edef\@p@sbbllx{\number\p@intvaluex}
                        \edef\@p@sbblly{\number\p@intvaluey}

                        \rotate@{\@p@sbburx}{\@p@sbbury}
                        \advance\p@intvaluex by \x@ffsetvalue
                        \advance\p@intvaluey by \y@ffsetvalue
                        \edef\@p@sbburx{\number\p@intvaluex}
                        \edef\@p@sbbury{\number\p@intvaluey}
%               swap LL and UR if necessary
%\typeout{rotated BB: \@p@sbbllx, \@p@sbblly, \@p@sbburx, \@p@sbbury}
                        {
                         \count0=\@p@sbbllx \count1=\@p@sbblly
                         \count2=\@p@sbburx \count3=\@p@sbbury
                         \dimen0=\@p@sbbllx sp\dimen1=\@p@sbblly sp
                         \dimen2=\@p@sbburx sp\dimen3=\@p@sbbury sp
                         \dimen203=\dimen2 \advance\dimen203 by -\dimen0
                         \dimen204=\dimen3 \advance\dimen204 by -\dimen1
                         \ifdim\dimen203<0sp 
                              \count203=\count2 \count2=\count0 
                              \count0=\count203 
                              \global\edef\@p@sbbllx{\number\count0}
                              \global\edef\@p@sbburx{\number\count2}
                         \fi
                         \ifdim\dimen204<0sp 
                               \count204=\count3
                               \count3=\count1
                               \count1=\count204
                               \global\edef\@p@sbblly{\number\count1}
                               \global\edef\@p@sbbury{\number\count3}
                         \fi
                        }
%\typeout{after swap BB: \@p@sbbllx, \@p@sbblly, \@p@sbburx, \@p@sbbury}
                \fi
                \count203=\@p@sbburx
                \count204=\@p@sbbury
                \advance\count203 by -\@p@sbbllx
                \advance\count204 by -\@p@sbblly
                \edef\@bbw{\number\count203}
                \edef\@bbh{\number\count204}
                %\typeout{ bbh = \@bbh, bbw = \@bbw }
}
%
%%% \in@hundreds performs #1 * (#2 / #3) correct to the hundreds,
%       then leaves the result in @result
%
\def\in@hundreds#1#2#3{\count240=#2 \count241=#3
                     \count100=\count240        % 100 is first digit #2/#3
                     \divide\count100 by \count241
                     \count101=\count100
                     \multiply\count101 by \count241
                     \advance\count240 by -\count101
                     \multiply\count240 by 10
                     \count101=\count240        %101 is second digit of #2/#3
                     \divide\count101 by \count241
                     \count102=\count101
                     \multiply\count102 by \count241
                     \advance\count240 by -\count102
                     \multiply\count240 by 10
                     \count102=\count240        % 102 is the third digit
                     \divide\count102 by \count241
                     \count200=#1\count205=0
                     \count201=\count200
                        \multiply\count201 by \count100
                        \advance\count205 by \count201
                     \count201=\count200
                        \divide\count201 by 10
                        \multiply\count201 by \count101
                        \advance\count205 by \count201
                     \count201=\count200
                        \divide\count201 by 100
                        \multiply\count201 by \count102
                        \advance\count205 by \count201
                     \edef\@result{\number\count205}
}
\def\compute@wfromh{
                % computing : width = height * (bbw / bbh)
                \in@hundreds{\@p@sheight}{\@bbw}{\@bbh}
                %\typeout{ \@p@sheight * \@bbw / \@bbh, = \@result }
                \edef\@p@swidth{\@result}
                %\typeout{w from h: width is \@p@swidth}
}
\def\compute@hfromw{
                % computing : height = width * (bbh / bbw)
                \in@hundreds{\@p@swidth}{\@bbh}{\@bbw}
                %\typeout{ \@p@swidth * \@bbh / \@bbw = \@result }
                \edef\@p@sheight{\@result}
                %\typeout{h from w : height is \@p@sheight}
}
\def\compute@handw{
                \if@height 
                        \if@width
                        \else
                                \compute@wfromh
                        \fi
                \else 
                        \if@width
                                \compute@hfromw
                        \else
                                \edef\@p@sheight{\@bbh}
                                \edef\@p@swidth{\@bbw}
                        \fi
                \fi
}
\def\compute@resv{
                \if@rheight \else \edef\@p@srheight{\@p@sheight} \fi
                \if@rwidth \else \edef\@p@srwidth{\@p@swidth} \fi
                %\typeout{rheight = \@p@srheight, rwidth = \@p@srwidth}
}
%               
%%% Compute any missing values
\def\compute@sizes{
        \compute@bb
        \compute@handw
        \compute@resv
}
%
%%% \epsfig
%%% usage : \epsfig{file=, height=, width=, bbllx=, bblly=, bburx=, bbury=,
%                       rheight=, rwidth=, angle=, clip=}
%
%%% "clip=" is a switch and takes no value, but the `=' must be present.
%------------------------------------------------------------------
%%% by the way, possible parameters to the PSfile= command in dvips are:
%%%                    llx
%%%                    lly
%%%                    urx
%%%                    ury
%%%                    rwi
%       hoffset The horizontal offset (default 0)
%       voffset The vertical offset (default 0)
%       hsize   The horizontal clipping size (default 612)
%       vsize   The vertical clipping size (default 792)
%       hscale  The horizontal scaling factor (default 100)
%       vscale  The vertical scaling factor (default 100)
%       angle   The rotation (default 0)
%------------------------------------------------------------------
\def\psfig#1{\vbox {
        % do a zero width hard space so that a single
        % \epsfig in a centering enviornment will behave nicely
        %{\setbox0=\hbox{\ }\ \hskip-\wd0}
        %
        \ps@init@parms
        \parse@ps@parms{#1}
        \ifnum\@p@scost<\@psdraft
                \typeout{[\@p@sfilefinal]}
                \if@verbose
                        \typeout{epsfig: using PSFIG macros}
                \fi
                \psfig@method
        \else
                \epsfig@draft
        \fi
}}

\def\epsfig#1{\vbox {
        % do a zero width hard space so that a single
        % \epsfig in a centering enviornment will behave nicely
        %{\setbox0=\hbox{\ }\ \hskip-\wd0}
        %
        \ps@init@parms
        \parse@ps@parms{#1}
        \ifnum\@p@scost<\@psdraft
                \typeout{[\@p@sfilefinal]}
                \if@clip\use@psfigtrue\fi
                \if@angle\use@psfigtrue\fi
                \ifuse@psfig
                        \if@verbose
                                \typeout{epsfig: using PSFIG macros}
                        \fi
                        \psfig@method
                \else
                        \if@verbose
                                \typeout{epsfig: using EPSF macros}
                        \fi
                        \epsf@method
                \fi
        \else
                \epsfig@draft
        \fi
}}

\def\epsf@method{%
        \epsfgetbb{\@p@sfile}%
        \epsfsetgraph{\@p@sfilefinal}
}
\def\psfig@method{%
        \compute@sizes
        \special{ps::[begin]  \@p@swidth \space \@p@sheight \space%
        \@p@sbbllx \space \@p@sbblly \space%
        \@p@sbburx \space \@p@sbbury \space%
        startTexFig \space }%
        \if@angle
                \special {ps:: \@p@sangle \space rotate \space} 
        \fi
        \if@clip
                \if@verbose
                        \typeout{(clipped to BB) }
                \fi
                \special{ps:: doclip \space }%
        \fi
        \special{ps: plotfile \@p@sfilefinal \space }%
        \special{ps::[end] endTexFig \space }%
        % Create the vbox to reserve the space for the figure%
        \vbox to \@p@srheight true sp{\hbox to \@p@srwidth true sp{\hss}\vss}
}
%
% draft figure, just reserve the space and print the
% path name.
\def\epsfig@draft{
\compute@sizes
\if@draftbox
        % Verbose draft: print file name in box
        % NOTE: fbox is a LaTeX command!
        \hbox{\fbox{\vbox to \@p@srheight true sp{
        \vss\hbox to \@p@srwidth true sp{ \hss \@p@sfilefinal \hss }\vss
        }}}
\else
        % Non-verbose draft
        \vbox to \@p@srheight true sp{%
        \vss\hbox to \@p@srwidth true sp{\hss}\vss}
\fi     
}
\catcode`\@=\atcode % return to previous catcode

\newcount\dcoef
\newcount\mcoef

\def\dessin(#1)(#2)(#3){{\epsfig{file=#1,width=#2truecm,height=#3truecm}}}

\mcoef=1
\dcoef=2

\gauche={S. Burckel}
\droite={The Parallel-Sequential Duality : Matrices and Graphs.}

\def\sq{{\downarrow}}

\centerline{\XIV The Parallel-Sequential Duality : Matrices and Graphs.}
\sk

\bn {\bf Abstract. Usually, mathematical objects have highly parallel interpretations.
In this paper, we consider them as sequential constructors of other objects.
In particular, we prove that every reflexive directed graph can be interpreted as a program that builds another and is
itself builded by another. That leads to some optimal memory computations, codings similar to modular decompositions and other
strange dynamical phenomenons.}

\bn
\ttparag{1. Introduction.}
\bn
In this paper, we will deal with matrices and finite directed graphs $G=(V,A)$ where $V$ is a finite list $(x_1,x_2,\dots,x_n)$ of vertices
and $A$ is a subset of $V\times V$ of arcs.
\bn
Usually, in classical mathematics, matrices are highly parallel objects. For computing the transformation of a vector $X$ by
a linear mapping $X:=E(X)$, a mathematician first keeps the initial vector $X$ in mind or in a safe place, and then computes
the successive images of each component of $X$ by the projections of $E$. Equivalently, he can also compute $Y:=E(X)$ and then $X:=Y$. In both cases,
one seems to have to build a copy of the initial vector $X$ in order to preserve its initial values for the whole computation that
will modify these values. Hence, the data size for this kind of computation processes is twice the size of the input data $X$.
A motivation of this paper is to show that one can compute the transformation $X:=E(X)$ without
any copy. In [1], we proved a similar result for the computation of boolean mappings.
Here, we will interpret matrices as sequential programs that perform sequences of assignments on the components of the initial vector.
In general, the mapping computed by this way is not the usual linear mapping represented. Now, a matrix admits a usual parallel interpretation and
a sequential interpretation. We are going to investigate this duality.
\bn

\Def{(matrices).}{ Let $K$ be a field. Given a square matrix $M$ of $M_{n,n}(K)$, $M_i$ is the $i$-th row vector of $M$ and $M_{i,j}$ is the $j$-th component
of the row $M_i$. The matrix $M$ is said {\it regular} if it only has $1$s on its diagonal, i.e., $M_{i,i}=1$ for every $i$.
A matrix $M'$ of $M_{n,n}(K)$ is {\it similar} to $M$ if it has the same values than $M$ excepted possibly on the diagonal, i.e., $M'_{i,j}\not=M_{i,j}\imp i=j$.
}

First, we recall the classical interpretation of a matrix in linear algebra.

\Def{(parallel mapping). }{ Let $K$ be a field. Every square matrix $M$ of $M_{n,n}(K)$ represents a linear mapping $M^\#$ from $K^n$ to $K^n$
that we will call the {\it parallel mapping of $M$} and defined by the following.
For every $X=(x_1,\dots,x_n)\in K^n$, the {\it parallel image of $X$ by $M$}
is the vector
$$M^\#(X)=(M_1.X,M_2.X,\dots,M_n.X)$$
This mapping $M^\#$ can be computed by the following straight-line program~:
\bn
{\tt
\sk Input : $X=(x_1,\dots,x_n)$
\sk
\sk For $i$ from $1$ to $n$ do $$y_i:=\sum_{j=1}^n M_{i,j}.x_j$$
\sk
\sk Output : $Y=(y_1,\dots,y_n)$
\sk
\sk
}
}

\bn
\ttparag{2. Sequentializing.}
\bn

Now, we define a dual interpretation of a matrix in linear algebra such that the image of a vector $X$ can be computed via
a sequence of linear transformations of this vector $X$ and only using this vector for the whole computation.
We obtain an "in situ" straight-line program.

\Def{(sequential mapping). }{ Let $K$ be a field. Every square matrix $M$ of $M_{n,n}(K)$ represents a linear mapping $M^\sq$ from $K^n$ to $K^n$
called the {\it sequential mapping of $M$} defined by the following.
For every $X=(x_1,\dots,x_n)\in K^n$, the {\it sequential image of $X$ by $M$}
is the resulting vector computed by the following straight-line program denoted $M^\pi$~:
\bn
{\tt
\sk Input : $X=(x_1,\dots,x_n)$
\sk
\sk For $i$ from $1$ to $n$ do $$x_i:=\sum_{j=1}^n M_{i,j}.x_j$$
\sk
\sk Output : $X=(x_1,\dots,x_n)$
\sk
\sk
}
Denote $M^s$ the matrix such that $(M^s)^\#=M^\sq$.
}

\bn
As we said above, one only uses the input vector $X$ for the computation of the sequential image, whereas in the
standard parallel interpretation, one uses a second vector $Y$ (in order to preserve the values of the input vector $X$).
\bn
For example, for $K=\RR$, the sequential mapping $M^\sq$ of the matrix
$$M=\left[\matrix{
0&1&2\cr
3&4&5\cr
6&7&8\cr
}\right]$$
maps every vector $(a,b,c)$ to $(b+2c,7b+11c,55b+97c)$ since the corresponding sequential program $M^\pi$ is~:

$$\matrix{
a&:=0.a+1.b+2.c&(=b+2c)&\cr
b&:=3.a+4.b+5.c&(=3(b+2c)+(4b)+(5c)=7b+11c)\cr
c&:=6.a+7.b+8.c&(=6(b+2c)+7(7b+11c)+8c=55b+97c)\cr
}$$

and one has~:

$$M^s=\left[\matrix{
0&1&2\cr
0&7&11\cr
0&55&97\cr
}\right]$$

\bn The sequential interpretation of a matrix leads to a natural notion in the context of directed graphs.

\Def{(graph construction). }{ Let $G=(V,A)$ be a finite directed graph where $V$ is a finite list $(x_1,x_2,\dots,x_n)$ of vertices
and $A$ is a subset of $V\times V$ of arcs. We say that $G$ {\it sequentially constructs} a directed graph $G'$ when
their respective adjacency matrices $M,M'$ in $M_{n,n}(F_2)$ satisfy : $$M^s=M'$$
}
\bn
That is a kind of decomposition of a graph. Instead of working with a graph $G$, one can consider a sequential constructor of $G$ or
the graph that $G$ constructs. We will see in the sequel some advantages. We are going to study the relations between these graphs and first of all,
their existence. Given a graph $G$, is there another graph $G'$ which is a sequential constructor of $G$ ? We will see that the answer is NO. However,
if one assumes the graph $G$ to be reflexive, the answer is YES.
\bn
In the formalism of matrices, observe that by definition, the parallel interpretation of $M^s$ is equal to the sequential interpretation of $M$.
For the other direction, a natural question is :

\sk{\it given a matrix $M$, is there a matrix $P$ with a sequential interpretation
equal to the parallel interpretation of $M$ ?}
\bn
Unfortunately, the answer is NO in general. A minimal example for $K=F_2$ is~:
$$M=\left[\matrix{
0&0\cr
1&0\cr
}\right]$$
If there were a matrix $P$ such that $P^\sq=M^\#$, since $P$ and $M$ have necessarily the same first rows (i.e., $P_1=M_1$),
the matrix $P$ must be some
$$\left[\matrix{
0&0\cr
k_1&k_2\cr
}\right]$$
and the sequential program $P^\pi$~:
{\tt
\sk $a:=0$
\sk $b:=k_1.a+k_2.b$
}
\sk
should transform every vector $(a,b)$ in $(0,a)$~: that is not possible.
\bn

However, we are not very far from a positive answer with the following preliminary result.

\Thm{}{Let $K$ be a field. For every linear mapping $E$ from $K^n$ to $K^n$ with $n>0$, the assignment $X:=E(X)$ is computed by a
straight-line program $P$ made of at most $2n-1$ linear assignments of the components of $X$. Moreover, the $n$ first steps of $P$
form a program $M^\pi$ for a matrix $M$ in $M_{n,n}(K)$.}

\Pr
We construct the code of the program $P$ in two parts.
\bn
Part 1.
{\tt
\sk Let $M$ be the matrix such that $M^\#=E$ and let $P$ be the empty program.
\sk
For $i$ from $1$ to $n$ do
\sk
\sk case 1. If $M_{ii}\not=0$
\sk Let $r_i:=0$
\sk Add to $P$ the instruction $[x_i:=M_i.X]$
\sk For every $k>i$ with $M_{ki}\not=0$ do $$M_k:=M_k+M_{ki}.M_{ii}^{-1}.(x_i-M_i)$$
\sk
\sk case 2. If $M_{ii}=0$ and for every $j>i$, $M_{ji}=0$
\sk Let $r_i:=0$
\sk Add to $P$ the instruction $[x_i:=M_i.X]$
\sk
\sk case 3. If $M_{ii}=0$ and there exists some $j>i$ with $M_{ji}\not=0$
\sk Let $r_i:=j$ for such a chosen $j$.
\sk Add to $P$ the instruction $[x_i:=(M_i-M_j).X]$
\sk For every $k>i$ with $M_{ki}\not=0$ do $$M_k:=M_k+M_{ki}.(-M_{ji})^{-1}.(x_i-M_i+M_j)$$
}
\bn
Part 2.
\bn
{\tt
for $i$ from $n-1$ downto $1$ do
\sk If $r_i=j>i$ then add to $P$ the instruction $[x_i:=x_i+x_j]$.
}
\bn
We prove the correction of this method, that is to say, that the program $P$ produced computes the transformation $X:=E(X)$ for every $X\in K^N$.
\bn
In Part 1, assume that at each step $i=1,\dots,n$ the rows $M_i,\dots,M_n$ represent the
respective images by $E$ of the components $x_i,\dots,x_N$ according to their current values at step $i$.
That is true for $i=1$ because $M^\#=E$.
At step $i$, the program makes an assignment of $x_i$.
\sk In case 1, $M_{ii}\not=0$ and $x_i$ receives the new value $a.x_i+b$ where $a=M_{ii}\not=0$. Hence, the initial value
of $x_i$ is not lost and was $a^{-1}.(x_i-b)$ where $x_i$ is its new value. In order to preserve the hypothesis on $M$,
one has to replace the references to the initial $x_i$ by this new linear function $a^{-1}.(x_i-b)$
in every $M_k$ for $k>i$. Hence every $M_k=M_{ki}x_i+F$ with $k>i$ and $M_{ki}\not=0$ should become
$$\eqalign{
M_{ki}.a^{-1}.(x_i-b)+F &=M_{ki}.a^{-1}.(x_i-M_i+a.x_i)+F\cr
				&=M_{ki}.a^{-1}.(x_i-M_i)+M_{ki}.x_i+F\cr
				&=M_{ki}.M_{ii}^{-1}.(x_i-M_i)+M_k
}$$
\sk In case 2, $M_{ii}=0$ and $M_{ii}^{-1}$ does not exist. Hence, the previous adjustment is not possible. However, this operation
is not necessary since the old value of $x_i$ is used in no $M_j$ for $j>i$.
\sk Case 3 is more critical. The old value of $x_i$ will be used again in some $M_j$ with $j>i$ but a direct assignment of $x_i$ would make
impossible to recover it because $M_{ii}=0$. The adopted solution here is to first transform $M_i$ in $M'_i=M_i-M_j$ in order to fail in case 1 since
$M'_{ii}=-M_{ji}\not=0$. The correction will be made in Part 2.
\bn
In Part 2, assume that for each $i=n-1,\dots,1$, the components $x_{i+1},\dots,x_n$ have received their correct values.
That is true for $i=n-1$, since the assignment of $x_n$ in Part 1 cannot fail in case 3 because $i=n$ is maximal (there is no $j>i$).
Hence $x_n$ has received its correct value by hypothesis in part 1.
If for $i$ such that $r_i=j>i$, by hypothesis, $x_j$ has received its correct value $E_j$. Hence,
the assignment $x_i:=x_i+x_j$ will perform $x_i:=x_i+E_j=E_i-E_j+E_j=E_i$ and $x_i$ also receives its correct value.
\Rp

Observe that the Part 1 of the produced program $P$ is exactly the sequential program $M^\pi$ and in Part 2 all the assignments are quite simple
and have the form $$x_i:=x_i+x_j$$
Hence, one can code Part 1 with the matrix $M$ and Part 2 with the list of numbers $[r_1,r_2,\dots,r_n]$.
\bn
For example, let us compute the program for the linear mapping $E:\RR^3\mapsto \RR^3$ with $E(a,b,c)=(c,a+b+c,3a+3b+2c)$

\sk For $i=1$, one has $M_{11}=0$ and we fall in case 3 since $M_{21}\not=0$ (and also $M_{31}\not=0$).
One chooses for instance $j=2$ and remember that for Part 2 with $r_1:=2$. Then performs $M_1=M_1-M_2=[-1,-1,0]$. The first assignment of the program is $a:=-a-b$.
Then $M_2$ becomes $[-1,0,+1]$ and $M_3$ becomes $[-3,0,2]$.

\sk For $i=2$, $M_{22}=0$ and we fall in case 2 since $M_{32}=0$ too. The second assignment of the program is $b:=-a+c$.
And $M_3$ remains equal to $[-3,0,2]$.

\sk For $i=3$, we are in case 1 and the third assignment of the program is $c:=-3a+2c$.

\sk In part 2, the only assignment added is $a:=a+b$. Finally, the computation program $P$ is

$$\eqalign{
a&:=-a-b\cr
b&:=-a+c\cr
c&:=-3a+2c\cr
a&:=a+b\cr
}$$

One can verify that $P$ computes $X:=E(X)$.

\sk Begin with $X=(a,b,c)$.

$$\matrix{
P		&|&	&X	&		\cr
		&|&a	&b	&c		\cr
a:=-a-b	&|&-a-b	&	&		\cr
b:=-a+c	&|&	&a+b+c&		\cr
c:=-3a+2c	&|&	&     &3a+3b+2c	\cr
a:=a+b	&|&c	&     &		\cr
}$$

\bn
From every $X=(a,b,c)\in\RR^3$, the program ends with $X=(c,a+b+c,3a+3b+2c)$ as expected.
\bn

The only problem for sequentializing is Case 3 that forces an adaptation in Part 2. Hence, one needs a larger structure than matrices
themselves in order to code this sequentializing~: for instance, a matrix AND the list of numbers $r_i$. In the previous example,
a possible coding of the sequentializing of

$$M=\left[\matrix{
0&0&1\cr
1&1&1\cr
3&3&2}\right]$$

is~:

$$\left[\matrix{
-1&-1&0\cr
-1&0&1\cr
-3&0&2}\right] [2,0,0]$$

For another example, the sequentializing of
$$M=\left[\matrix{
0&0\cr
1&0\cr
}\right]$$
leads to the following straight-line program~:
$$\eqalign{
a&:=-a\cr
b&:=-a\cr
a&:=a+b
}$$
and a possible coding is~:
$$\left[\matrix{
-1&0\cr
-1&0\cr
}\right] [2,0]$$

\bn
Now we consider matrices up to reflexivity (i.e., up to their diagonals).
With this restriction, the previous sequentializing method can be simplified such that only Cases 1 occur.
Hence, the sequentializing program is only made of $n$ steps where $n$ is the dimension of
the given matrix (or the graph).

\Thm{}{Let $K$ be a field. Let $M$ be a matrix in $M_{n,n}(K)$.
There exists a regular matrix ${}^dM$ in $M_{n,n}(K)$ such that ${}^dM^s$ is similar to $M$.
}

\Pr
The construction is similar than in the previous proof. However, the adjustments consist in the replacement at step $i$ of
a critical case $M_{i,i}=0$ by forcing $M_{i,i}:=1$. We explicit now this procedure in the particular case $K=F_2$.
\bn
{\tt
Begin with ${}^dM=M$ and $P$ is the empty program.
\sk For $i$ from $1$ to $N$ do
\sk Set ${}^dM_{i,i}:=0$.
\sk Add to $P$ the instruction $[x_i:=x_i+{}^dM_i.X]$
\sk For every $k>i$ with ${}^dM_{ki}=1$ do $${}^dM_k:={}^dM_k+{}^dM_i$$
\sk Set ${}^dM_{i,i}:=1$.
\sk
\sk
}
As previously, the construction implies that at each step $i$, the rows ${}^dM_i,\dots,{}^dM_n$ represent the images (up to reflexivity) of the
initial vector according to its current values.
\Rp
\bn
For example, for
$$M=\left[\matrix{
0&1&1\cr
1&1&0\cr
1&0&1}\right]$$
one obtains~:
\sk for $i=1$, one adds to $P$ the instruction $[a:=a+b+c]$ and ${}^dM$ becomes
$${}^dM=\left[\matrix{
1&1&1\cr
1&0&1\cr
1&1&0}\right]$$
\sk for $i=2$, one adds to $P$ the instruction $[b:=a+b+c]$ and ${}^dM$ becomes
$${}^dM=\left[\matrix{
1&1&1\cr
1&1&1\cr
0&1&1}\right]$$
\sk for $i=3$, one adds to $P$ the instruction $[c:=b+c]$
\sk We have obtained the regular matrix
$${}^dM=\left[\matrix{
1&1&1\cr
1&1&1\cr
0&1&1}\right]$$
where $${}^dM^s=\left[\matrix{
1&1&1\cr
1&0&0\cr
1&0&1}\right]$$
is similar to $M$.
\bn
\ttparag{3. Dynamics.}
\bn

Some strange dynamical systems can be deduced from the above. Let $G$ be a reflexive directed graph. Hence its adjacency matrix $M$
is boolean and regular.
Let us construct the sequence of matrices $M=M_0,M_1,M_2,\dots$ where $M_{i+1}={}^dM_i$.
By construction, these matrices are all boolean and regular. Since the set containing these matrices is finite,
there must exists some cycle, i.e. some pair $0\le p<q$ with $M_p=M_q$ and where $p$ is taken minimal.
\sk We claim that $p=0$.
\sk Observe that this sequence of matrices can be constructed in the other direction. From $M_{i+1}$, build the matrix $M_{i+1}^s$. With the previous
result, $M_{i+1}^s={}^dM_i^s$ is similar to $M_i$. But $M_i$ is regular. Hence, by putting only $1$s on the diagonal of $M_{i+1}^s$, one obtains a deterministic way
to transform $M_{i+1}$ to $M_i$. Denote $\phi$ this mapping on regular matrices that takes a matrix $M$, computes $M^s$ and puts $1$s on its diagonal.
\sk
Now assume that $p>0$. The matrix $M_p$ would satisfy $\phi(M_p)=M_{p-1}$ and $\phi(M_p)=M_q$. Since $\phi$ is a mapping, one obtains $M_{p-1}=M_q$~: a contradiction
with the minimality of $p$.
\bn
For example, the maximal cycle on matrices of $M_{4,4}(F_2)$ has length $18$ and from the matrix~:

$$M_0=\left[\matrix{
1& 1& 1& 1& 1& 0& 1& 0& 1& 1\cr
1& 1& 1& 0& 1& 1& 0& 0& 1& 0\cr
0& 1& 1& 0& 1& 0& 0& 1& 0& 1\cr
0& 0& 1& 1& 1& 1& 1& 0& 0& 1\cr
0& 0& 0& 1& 1& 1& 1& 0& 1& 1\cr
1& 1& 0& 0& 0& 1& 0& 0& 1& 1\cr
0& 1& 0& 1& 1& 0& 1& 0& 0& 1\cr
1& 1& 0& 1& 0& 0& 1& 1& 0& 0\cr
1& 1& 0& 0& 1& 1& 1& 1& 1& 1\cr
1& 0& 0& 1& 1& 1& 0& 0& 0& 1\cr
}\right]$$

one obtains a cycle of length $13122$.

\bn
\ttparag{4. Relations.}
\bn
We have seen that there exist matrices $M$ that do not admit some matrice $P$ so that $P^s=M$. However,
every matrix $M$ has a unique sequential mapping $M^\sq$. Hence there exist pairs of matrices $(P,P')$ that construct the same matrix. Let us
investigate in this paragraph these relations.

\Def{(equivalence).}{ Two matrices $M,W$ are {\it sequentially equivalent} if their sequential mappings $M^\sq$ and $W^\sq$ are equal.
}

For example, the two matrices of $M_{4,4}(\RR)$~:

$$M=\left[\matrix{
1&2&3&4\cr
1&0&0&0\cr
0&1&0&0\cr
0&0&1&0\cr
}\right]\ \ W=\left[\matrix{
1&2&3&4\cr
1&0&0&0\cr
1&0&0&0\cr
1&0&0&0\cr
}\right]
$$
are sequentially equivalent since the two mappings $M^\sq$ and $W^\sq$ are equal to
$$(a,b,c,d)\mapsto(a+2b+3c+4d,a+2b+3c+4d,a+2b+3c+4d,a+2b+3c+4d)$$
However, their sequential programs are different~:
$$\matrix{
M^\pi&W^\pi\cr
a:=a+2b+3c+4d&a:=a+2b+3c+4d\cr
b:=a&b:=a\cr
c:=b&c:=a\cr
d:=c&d:=a
}$$
\bn

Of course, this notion can be exported in the context of directed graphs.

\Def{(equivalence).}{ Two directed graph $G,G'$ are {\it sequentially equivalent} if they sequentially construct the same graph.
In that case, we write $G\equiv G'$.
}

\Prop{(chains).}{ Every directed graph $G=(V=(x_1,\dots,x_n),A)$ that admits an induced subgraph $H$ on vertices $(x_p,x_{p+1},\dots,x_q)$ which is
a chain (i.e., $(x_{i+1},x_i)\in A$ for $p\le i<q$), is sequentially equivalent to the graph $G'$ where an arc $(x_{i+1},x_i)$ of this chain is replaced
by an arc $(x_{i+1},x_j)$ with $p\le j<i+1$.}

\Pr
The corresponding sequential program performs~:
{\tt
\sk $x_{p+1}:=x_p$
\sk $x_{p+2}:=x_{p+1}$
\sk $x_{p+3}:=x_{p+2}$
\sk\dots
\sk $x_{q}:=x_{q-1}$
}
\sk and computes the same mapping than the program~:
{\tt
\sk $x_{p+1}:=x_p$
\sk $x_{p+2}:=x_p$
\sk $x_{p+3}:=x_p$
\sk\dots
\sk $x_{q}:=x_p$
}
\Rp

For example~:
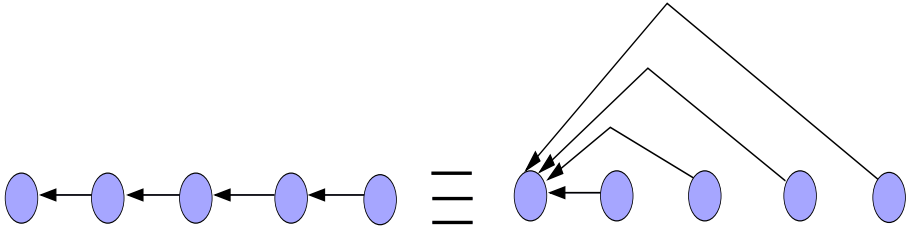
$$\dessin(eq1.eps)(12)(3)$$

The complexity of graphs can be considerably reduced up to sequential equivalence. For instance, when a graph contains a maximal directed acyclic subgraph.

\Prop{(linear orders).}{ Every directed graph $G=(V=(x_1,\dots,x_n),A)$ that admits an induced subgraph $H$ on vertices $(x_p,x_{p+1},\dots,x_q)$ which is
a linear order (i.e., $(x_j,x_i)\in A$ for $p\le i<j\le q$), is sequentially equivalent to the graph $G'$ where one removes all the arcs of $H$ excepted the arc
$(x_{p+1},x_p)$.}

\Pr
The corresponding sequential program performs~:
{\tt
\sk $x_{p+1}:=x_p$
\sk $x_{p+2}:=x_{p+1}+x_p$
\sk $x_{p+3}:=x_{p+2}+x_{p+1}+x_p$
\sk\dots
\sk $x_{q}:=x_{q-1}+\dots+x_{p+2}+x_{p+1}+x_p$
}
\sk and computes the same mapping than the program~:
{\tt
\sk $x_{p+1}:=x_p$
\sk $x_{p+2}:=0$
\sk $x_{p+3}:=0$
\sk\dots
\sk $x_{q}:=0$
}
\Rp

For example~:
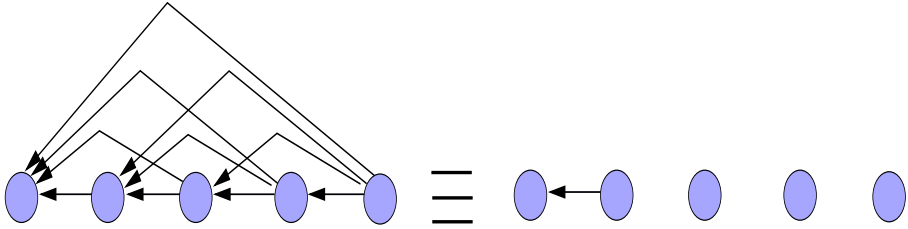
$$\dessin(eq2.eps)(12)(3)$$

\bn

\ttparag{5. Conclusion.}

The sequential interpretation of matrices we presented here is a powerful tool for computations of linear transformations $X:=E(X)$
using minimal memory~: one only uses the input data memory $X$ for that computation. In the case of graphs, and more generally directed graphs,
the sequential interpretation defines a graph decomposition~: roughly speaking, every vertex of rank $i$ can "use" the connecting work of vertices of ranks
$j<i$. Moreover, sequential equivalences in graphs enable some simplifications in the codings. A complete axiomatization of these relations could be done.
\bn
Under a dynamical system point of view, iterating sequential interpretations
leads to strange phenomena which should also merit particular investigations.
\bn
Let us mention that other choices of construction are possible. Here, the sequential constructor ${}^dM$ of a matrix $M$ is chosen to have only $1s$ on
its diagonal. However, every sequence of invertible elements in the field $K$ is possible. Moreover, similar constructions can be done on finite sets in
general. We have seen here two sequential decompositions of a linear mapping on $K^n$. The first one consists in $2n-1$ steps where the last $n-1$ steps
are quite minimal and coded by assignments $x_i:=x_i+x_j$ with $j>i$ or by a sequence of integers $[r_1,r_2,\dots,r_{n-1}]$.
The second one just consists in $n$ steps but the diagonal of the matrix can be modified. Let us propose a third possibility, still in $n$ steps, where each row of the matrix is completely preserved.
However, the order of the rows may change~: during the process, when $M_{i,i}=0$ and there exists $j>i$ with $M_{j,i}\not=0$, then exchange the two rows $M_i$ and $M_j$
and proceed like in the two other methods by substitution.
Now, the coding of the result is a new matrix and a permutation $\sigma$ of $S_n$. This method preserves the informative structure of the mapping.
The similar operation on directed graphs, up to a permutation of vertices, preserves the numbers $d^-(x)$ of incoming arcs to each vertex $x$.
More precisely, there exists a permutation of vertices $\sigma$ such that every arc $(x,y)$ is transformed in an arc $(\sigma(x),y)$.
The semantic of this operation also merits other investigations.

\bn
\ttparag{References.}
%\Ref  [abrev]; auteur; titre; reste;
%\Reff [abrev]; auteur; titre; journal; numero; annee; page;

\def \it#1{#1}

\parindent=1.8truecm \def\Ref#1; #2; #3; #4;{ \sgni\item{\hbox to 1.6truecm{ [#1]\hfill}}{\sc #2},
{\sl #3}, #4. } \def\Reff#1; #2; #3; #4; #5; #6; #7;{ \sgni\item{\hbox to 1.6truecm{
[#1]\hfill}}{\sc #2}, {\sl #3}, #4 {\bf #5} (#6) #7. }

\Reff {1}; S. Burckel; Closed Iterative Calculus; Theoretical Computer Science; 158; 1996; 371--378;

\bn

\bn

\sk\ni Serge Burckel.
\sk\ni INRIA-LORIA,
\sk\ni Campus Scientifique
\sk\ni BP 239
\sk\ni 54506 Vandoeuvre-l\`es-Nancy Cedex
\sk\ni sergeburckel@orange.fr
\bn

{\it This work has been partially supported by the CNRS-ANR project "graphs decompositions and algorithms".}

\end